\documentclass[a4paper,12pt]{amsart}
\usepackage{amssymb}
\usepackage{ifthen}
 \usepackage[dvips]{graphicx}
\nonstopmode \numberwithin{equation}{section}
\usepackage{amssymb}
\usepackage{ifthen}
\usepackage{graphicx}
\usepackage{amsmath}
\usepackage[T1]{fontenc} 
\usepackage[utf8]{inputenc}
\usepackage[usenames,dvipsnames]{color}
\usepackage{color}
\usepackage[english]{babel}
\usepackage{fancyhdr}
\usepackage{fancybox}
\usepackage{tikz}

\nonstopmode \numberwithin{equation}{section}
\theoremstyle{plain}
\newtheorem{thm}[equation]{Theorem}
\newtheorem{cor}[equation]{Corollary}
\newtheorem{lem}[equation]{Lemma}
\newtheorem{prop}{Proposition}

\newtheorem{conj}{Conjecture}

\theoremstyle{definition}
\newtheorem{defn}{Definition}[section]

\newtheorem{prob}{Problem}
\newtheorem{rem}{Remark}[section]

\newcounter{minutes}
\divide\time by 60
\newcounter{hours}
\multiply\time by 60
\addtocounter{minutes}{-\time}

\newcounter {own}
\def\theown {\thesection       .\arabic{own}}

\newenvironment{pf}[1][]{%
 \vskip 3mm
 \noindent
 \ifthenelse{\equal{#1}{}}%
  {{\slshape Proof. }}%
  {{\slshape #1.} }%
 }%
{\qed\bigskip}

\newcounter{alphabet}

\newcommand{\real}{{\operatorname{Re}\,}}

\def\be{\begin{equation}}
\def\ee{\end{equation}}

\newcommand{\bee}{\begin{enumerate}}
\newcommand{\eee}{\end{enumerate}}

\newcommand{\blem}{\begin{lem}}
\newcommand{\elem}{\end{lem}}
\newcommand{\bthm}{\begin{thm}}
\newcommand{\ethm}{\end{thm}}
\newcommand{\bcor}{\begin{cor}}
\newcommand{\ecor}{\end{cor}}
\newcommand{\beg}{\begin{examp}}
\newcommand{\eeg}{\end{examp}}
\newcommand{\begs}{\begin{examples}}
\newcommand{\eegs}{\end{examples}}

\newcommand{\bdefn}{\begin{defn}}
\newcommand{\edefn}{\end{defn}}

\newcommand{\bprob}{\begin{prob}}
\newcommand{\eprob}{\end{prob}}
\newcommand{\bei}{\begin{itemize}}
\newcommand{\eei}{\end{itemize}}

\newcommand{\bcon}{\begin{conj}}
\newcommand{\econ}{\end{conj}}
\newcommand{\bcons}{\begin{conjs}}
\newcommand{\econs}{\end{conjs}}
\newcommand{\bprop}{\begin{prop}}
\newcommand{\eprop}{\end{prop}}
\newcommand{\br}{\begin{rem}}
\newcommand{\er}{\end{rem}}
\newcommand{\brs}{\begin{rems}}
\newcommand{\ers}{\end{rems}}
\newcommand{\bo}{\begin{obser}}
\newcommand{\eo}{\end{obser}}
\newcommand{\bos}{\begin{obsers}}
\newcommand{\eos}{\end{obsers}}
\newcommand{\bpf}{\begin{pf}}
\newcommand{\epf}{\end{pf}}
\newcommand{\ba}{\begin{array}}
\newcommand{\ea}{\end{array}}
\newcommand{\beq}{\begin{eqnarray}}
\newcommand{\beqq}{\begin{eqnarray*}}
\newcommand{\eeq}{\end{eqnarray}}
\newcommand{\eeqq}{\end{eqnarray*}}

\begin{document}

\title{Improved Bohr inequalities for certain class of harmonic univalent functions}

\author{Molla Basir Ahamed}
\address{Molla Basir Ahamed,
	School of Basic Science,
	Indian Institute of Technology Bhubaneswar,
	Bhubaneswar-752050, Odisha, India.}
\email{mba15@iitbbs.ac.in}

\author{Vasudevarao Allu}
\address{Vasudevarao Allu,
School of Basic Science,
Indian Institute of Technology Bhubaneswar,
Bhubaneswar-752050, Odisha, India.}
\email{avrao@iitbbs.ac.in}

\author{Himadri Halder}
\address{Himadri Halder,
School of Basic Science,
Indian Institute of Technology Bhubaneswar,
Bhubaneswar-752050, Odisha, India.}
\email{hh11@iitbbs.ac.in}

\subjclass[{AMS} Subject Classification:]{Primary 30C45, 30C50, 30C80}

\def\thefootnote{}
\footnotetext{ {\tiny File:~\jobname.tex,
printed: \number\year-\number\month-\number\day,
          \thehours.\ifnum\theminutes<10{0}\fi\theminutes }
} \makeatletter\def\thefootnote{\@arabic\c@footnote}\makeatother

\begin{abstract}
Let $ \mathcal{H} $ be the class of complex-valued harmonic mappings $ f=h+\bar{g}$ defined in the unit disk $ \mathbb{D} : =\{z\in\mathbb{C} : |z|<1\} $, where $ h $ and $ g $ are analytic functions in $ \mathbb{D} $ with the normalization $ h(0)=0=h^{\prime}(0)-1 $ and $ g(0)=0 $. 
Let $ \mathcal{H}_{0}=\{f=h+\bar{g}\in\mathcal{H} : g^{\prime}(0)=0\}. $ Ghosh and Vasudevrao \cite{Ghosh-Vasudevarao-BAMS-2020} have studied the following interesting harmonic univalent class $ \mathcal{P}^{0}_{\mathcal{H}}(M) $ which is defined by
 $$\mathcal{P}^{0}_{\mathcal{H}}(M) :=\{f=h+\overline{g} \in \mathcal{H}_{0}: \real (zh^{\prime\prime}(z))> -M+|zg^{\prime\prime}(z)|,\; z \in \mathbb{D}\; \mbox{and}\;\; M>0\}.
 $$  In this paper, we obtain the sharp Bohr-Rogosinski inequality, improved Bohr inequality, refined Bohr inequality and Bohr-type inequality for the class $ \mathcal{P}_{\mathcal{H}}^{0}(M) $.
\end{abstract}

\maketitle
\pagestyle{myheadings}
\markboth{Molla Basir Ahamed, Vasudevarao Allu and  Himadri Halder}{Improved Bohr inequalities for certain class of harmonic univalent functions}

\section{Introduction}

The classical Bohr inequality (see \cite{Bohr-1914}) states that if $ f $ be an analytic function with the power series representation $ f(z)=\sum_{n=0}^{\infty}a_nz^n $ in $ \mathbb{D} $ where $ \mathbb{D}:=\{z\in\mathbb{C} : |z|<1\} $ such that $ |f(z)|\leq 1 $ for all $ z\in\mathbb{D} $, then 
\begin{equation}\label{e-1.1a}
	\sum_{n=0}^{\infty}|a_n|r^n\leq 1,\;\; \text{for all}\;\; |z|=r\leq\frac{1}{3}, 
\end{equation}
and the constant $ 1/3 $ can not be improved. The constant $ r_0=1/3 $ is known as Bohr's radius, while the inequality $ \sum_{n=0}^{\infty}|a_n|r^n\leq 1 $ is known as Bohr inequality. Bohr actually obtained the inequality \eqref{e-1.1a} for $ r\leq 1/6 $ and later Weiner, Riesz and Schur have independently proved it to $ 1/3 $. \par The Bohr inequality can be written  in terms of distance formulation as follows
\begin{equation*}
	d\left(\sum_{n=0}^{\infty}|a_nz^n|,|a_0|\right)=\sum_{n=1}^{\infty}|a_nz^n|\leq 1-f(0)=d(f(0),\partial\mathbb{D}),
\end{equation*}
where $ d $ is the Euclidean distance and $ \partial\mathbb{D} $ is the boundary of the unit disk $ \mathbb{D}. $ The notion of Bohr inequality can be generalized to any domain $ \Omega $ to find the largest radius $ r_{_{\Omega}}>0 $ such that 
\begin{equation}\label{e-1.1}
	d\left(\sum_{n=0}^{\infty}|a_nz^n|,|a_0|\right)=\sum_{n=1}^{\infty}|a_nz^n|\leq d(f(0),\partial\Omega)
\end{equation}
holds for all $ |z|=r\leq r_{_{\Omega}} $, and for all functions analytic in $ \mathbb{D} $ and such that $ f(\mathbb{D})\subseteq\Omega. $ Interestingly enough, it was exhibited in \cite{aizn-2007} that if $ \Omega $ is convex then the inequality \eqref{e-1.1} holds for all $ |z|\leq 1/3 $, and this radius is the best possible. Thus if $ \Omega $ is convex, then $ r_{\Omega} $ coincides with Bohr's radius and most notably, the radius does not depend on $ \Omega $. Abu-Muhana \cite{Abu-2010} established a result showing that when $ \Omega $ is a simply connected domain with $ f(\mathbb{D})\subset\Omega $, then the inequality \eqref{e-1.1} holds for $ |z|\leq 3-2\sqrt{2}=0,1715... $ and this radius is sharp for the Koebe function $ k(z)=z/(1-z)^2. $ 
\vspace{3mm}

Operator algebraists began to get interested in the inequality after Dixon \cite{Dixon & BLMS & 1995} exposed a connection between the inequality and the characterization of Banach algebras that satisfy von Neumann inequality. The generalization of Bohr's theorem for different classes of analytic functions becomes now a days an active research area (see \cite{kayumov-2018-b,kayumov-2018-c,Liu-2020}). For example, for the holomorphic functions, Aizenberg \textit{et al.} \cite{aizenberg-2001}, Aytuna and Djakov \cite{Ayt & Dja & BLMS & 2013} have studied Bohr phenomenon, and for the class of starlike logarithmic mappings, Ali \textit{et al.} \cite{Ali & Abdul & Ng & CVEE & 2016} found Bohr radius. In $ 2018 $, Ali and Ng \cite{Ali & Ng & CVEE & 2016} extended the classical Bohr inequality in Poincar$ \acute{e} $ disk model of the hyperbolic plane. In $ 2018 $, Kayumov and Ponnusamy \cite{kayumov-2018-b} introduced the notion of $ p $-Bohr radius for harmonic functions, and established result obtaining $ p $-Bohr radius for the class of odd analytic functions. Powered Bohr radius for the class of all self analytic maps on $ \mathbb{D} $ has been studied in \cite{Kay & Pon & AASFM & 2019} while several different improved versions of the classical Bohr inequality were proved in \cite{kayumov-2018-c}. In this connection, Kayumov \textit{et al.} \cite{kay & Pon & Sha & MN & 2018} ascertained Bohr radius for the class of analytic Bloch functions and also for certain harmonic Bloch functions.
\vspace{3mm}

The main aim of this paper is to establish several improved versions of Bohr inequality, Refined-Bohr inequality and Bohr-Rogosinski inequality, finding the corresponding sharp radius for the class $ \mathcal{P}^{0}_{\mathcal{H}}(M) $ which has been studied by Ghosh and Vasudevarao in \cite{Ghosh-Vasudevarao-BAMS-2020} 
$$\mathcal{P}^{0}_{\mathcal{H}}(M)=\{f=h+\overline{g} \in \mathcal{H}_{0}: \real (zh^{\prime\prime}(z))> -M+|zg^{\prime\prime}(z)|, \; z \in \mathbb{D}\; \mbox{and }\; M>0\},
$$ where \begin{equation}\label{e-1.2}
	f(z)=h(z)+\overline{g(z)}=z+\sum_{n=2}^{\infty}a_nz^n+\overline{\sum_{n=2}^{\infty}b_nz^n}.
\end{equation}

\noindent To study Bohr inequality and Bohr radius for functions in $ \mathcal{P}^{0}_{\mathcal{H}}(M) $, we require the coefficient bounds and growth estimate of functions in $ \mathcal{P}^{0}_{\mathcal{H}}(M) $. We have the following result on the coefficient bounds and and growth estimate for functions in $ \mathcal{P}^{0}_{\mathcal{H}}(M) $.
\begin{lem} \label{lem-1.2} \cite{Ghosh-Vasudevarao-BAMS-2020}
Let $f=h+\overline{g}\in \mathcal{P}^{0}_{\mathcal{H}}(M)$ for some $M>0$ and be of the form \eqref{e-1.2}. Then for $n\geq 2,$ 
\begin{enumerate}
\item[(i)] $\displaystyle |a_n| + |b_n|\leq \frac {2M}{n(n-1)}; $\\[2mm]

\item[(ii)] $\displaystyle ||a_n| - |b_n||\leq \frac {2M}{n(n-1)};$\\[2mm]

		\item[(iii)] $\displaystyle |a_n|\leq \frac {2M}{n(n-1)}.$
	\end{enumerate}
	The inequalities  are sharp with extremal function   $f$ given by 
	$f^{\prime}(z)=1-2M\, \ln\, (1-z) .$	
\end{lem}
\begin{lem}\cite{Ghosh-Vasudevarao-BAMS-2020}\label{lem-1.3}
	Let $f \in \mathcal{P}^{0}_{\mathcal{H}}(M)$. Then 
	\begin{equation} \label{e-1.4}
		|z| +2M \sum\limits_{n=2}^{\infty} \dfrac{(-1)^{n-1}|z|^{n}}{n(n-1)} \leq |f(z)| \leq |z| + 2M \sum\limits_{n=2}^{\infty} \dfrac{|z|^{n}}{n(n-1)}.
	\end{equation}
	Both  inequalities are sharp for the function $f_{M}$ given by $f_{M}(z)=z+ 2M \sum\limits_{n=2}^{\infty} \dfrac{z^n}{n(n-1)} .
	$
\end{lem}
The organization of this paper as follows: In section 2, we prove the sharp Bohr-Rogosinski radius for the class $ \mathcal{P}^{0}_{\mathcal{H}}(M) $.  In section 3, we prove the sharp results on improved-Bohr radius for the class $ \mathcal{P}^{0}_{\mathcal{H}}(M) $. In section 4, we prove the sharp results on refined Bohr inequality the sharp results on Bohr-type inequalities. In section 6, we give the proof of all the main results of this paper.
 \section{Bohr-Rogosinski radius for the class $ \mathcal{P}^{0}_{\mathcal{H}}(M) $}
 We prove the following sharp Bohr-Rogosinski inequality for function in the class $ \mathcal{P}^{0}_{\mathcal{H}}(M) $.
 \begin{thm}\label{th-2.1}
 	Let $ f\in \mathcal{P}^{0}_{\mathcal{H}}(M) $ be given \eqref{e-1.2}. Then for $ N\geq 2 $, we have
 	\begin{equation}\label{e-2.2}
 		|f(z)|+\sum_{n=N}^{\infty}(|a_n|+|b_n|)|z|^n\leq d\left(f(0),\partial f(\mathbb{D})\right)
 	\end{equation} 
 	for $ r\leq r_{_N}(M) $, where $ r_{_N}(M) $ is the smallest root of the equation
 	\begin{equation}\label{e-2.3}
 		r-1+2M\left(2r-1+2(1-r)\ln\,((1-r))-\sum_{n=2}^{N-1}\frac{r^n}{n(n-1)}+\ln 4\right)=0.
 	\end{equation}
 	The radius $ r_{_N}(M) $ is the best possible.
 \end{thm}
 
  \begin{rem} In the study of the roots $ r_{_N}(M) $ of the equation \eqref{e-2.3}, following interesting facts can be observed
 	\begin{enumerate}
 		\item[(i).] For the values of 
 		\begin{equation*}
 			M\geq\frac{1}{2(2\ln 2-1)}
 		\end{equation*} the equation 
 		\begin{equation*}
 			H_{_{N,M}}(r):=r-1+2M\left(2r-1+2(1-r)\ln\,((1-r))-\sum_{n=2}^{N-1}\frac{r^n}{n(n-1)}+\ln 4\right)=0
 		\end{equation*} does not have any solutions in $ (0,1) $.
 		\item[(ii).] For all odd $ N $, the equation $ H_{_{N,M}}(r)\neq 0 $ for $ r\in (0,1) $ although the function $ H_{_{N,M}}(r) $ has some finite global maximum at $ r=1 $.
 		\item[(iii).] For all even $ N $, the equation $ H_{_{N,M}}(r)=0 $ has no roots in $ \mathbb{R} $.
 	\end{enumerate}
 \end{rem}
 
 \begin{figure}[!htb]
 	\begin{center}
 		\includegraphics[width=0.60\linewidth]{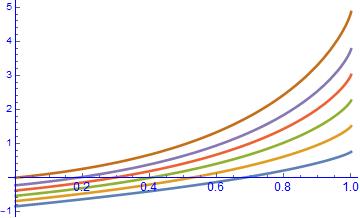}
 	\end{center}
 	\caption{The graph of $ r_3(M) $ of \eqref{e-2.2} when $ M<{1}/{2(2\ln 2-1)} $.}
 \end{figure} 
 \begin{table}[ht]
 	\centering
 	\begin{tabular}{|l|l|l|l|l|l|l|}
 		\hline
 		$\;\;M$& $\;\;0.2$&$\;\;0.4$& $\;\;0.6$& $\;\; 0.8 $& $\;\;1.0$&$\;\;1.29$\\
 		\hline
 		$r_3(M)$& $0.683 $&$0.527 $& $0.405$& $0.296$& $0.187 $&$0.003 $\\
 		\hline
 		$r_5(M)$& $0.702 $&$0.541 $& $0.414$& $0.301$& $0.189 $&$0.003 $\\
 		\hline
 		$r_{7}(M)$& $0.705 $&$0.542 $& $0.414$& $0.300$& $0.189 $&$0.003 $\\
 		\hline
 		$r_{9}(M)$& $0.706 $&$0.5429 $& $0.414$& $0.300$& $0.1891 $&$0.003 $\\
 		\hline
 	\end{tabular}
 	\vspace{2mm}
 	\caption{The roots $ r_{_N}(M) $ for $ N=3, 5, 7, 9 $ when \\$ M<{1}/{2(2\ln 2-1)} $.}
 	\label{tabel-1}
 \end{table}

 \begin{figure}[!htb]
 	\begin{center}
 		\includegraphics[width=0.60\linewidth]{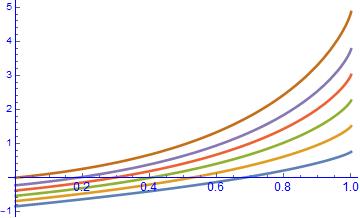}
 	\end{center}
 	\caption{The root $ r_{_8}(M) $ of \eqref{e-2.2} when $ M<{1}/{2(2\ln 2-1)} $.}
 \end{figure}
 
 By considering power of $ |f(z)| $ for the functions in the class $\mathcal{P}^{0}_{\mathcal{H}}(M) $, we obtain the following sharp result showing that the radius is different to compare with Theorem \ref{th-2.1}.
 \begin{thm}\label{th-2.4}
 	Let $ f\in \mathcal{P}^{0}_{\mathcal{H}}(M) $ be given by \eqref{e-1.2}. Then for $ N\geq 2 $,
 	\begin{equation}\label{e-2.5}
 		|f(z)|^2+\sum_{n=N}^{\infty}(|a_n|+|b_n|)|z|^n\leq d\left(f(0),\partial f(\mathbb{D})\right)
 	\end{equation} 
 	for $ r\leq r_{_N}(M) $, where $ r_{_N}(M) $ is the smallest root of the equation
 	\begin{align}
 		\label{e-2.6}&
 		(r+2M(r+(1-r)\ln(1-r)))^2-1\\&\nonumber+2M\left(r-1+(1-r)\ln(1-r)-\sum_{n=1}^{N-1}\frac{r^n}{n(n-1)}+\ln 4\right)=0.
 	\end{align}
 	The radius $ r_{_N}(M) $ is the best possible.
 \end{thm} 
 \begin{table}[ht]
	\centering
	\begin{tabular}{|l|l|l|l|l|l|l|}
		\hline
		$\;\;M$& $\;\;0.2$&$\;\;0.4$& $\;\;0.6$& $\;\; 0.8 $& $\;\;1.0$&$\;\;1.29$\\
		\hline
		$r_{_4}(M)$& $0.697 $&$0.537 $& $0.412$& $0.300$& $0.189 $&$0.003 $\\
		\hline
		$r_{_6}(M)$& $0.704 $&$0.542 $& $0.415$& $0.301$& $0.189 $&$0.003 $\\
		\hline
		$r_{_8}(M)$& $0.706 $&$0.542 $& $0.414$& $0.300$& $0.189 $&$0.003 $\\
		\hline
	\end{tabular}
	\vspace{1mm}
	\caption{The roots $ r_{_N}(M) $ for $ N=4, 6, 8 $  when $ M<{1}/{2(2\ln 2-1)} $.}
	\label{tabel-2}
\end{table} 
\begin{table}[ht]
	\centering
	\begin{tabular}{|l|l|l|l|l|l|l|l|l|l|}
	\hline
	$\;\;M$& $\;\;0.2$&$\;\;0.3$& $\;\;0.5$& $\;\; 0.6 $& $\;\;0.8$&$\;\;0.9$& $\;\;1.0$& $ \;\;1.1 $& $ \;\;1.29 $ \\
	\hline
	$r_{_3}(M)$& $0.739 $&$0.668 $& $0.559$& $0.512$& $0.425 $&$0.381 $& $0.334$& $0.279$ &$ 0.053 $\\
	\hline
	$r_{_5}(M)$& $0.751 $&$0.681 $& $0.517$& $0.524$& $0.435 $&$0.390 $& $0.342$& $0.286$ &$ 0.054 $\\
	\hline
	$r_{_6}(M)$& $0.753 $&$0.682 $& $0.572$& $0.525$& $0.436 $&$0.391 $& $0.342$& $0.286$ &$ 0.054 $\\
	\hline
	$r_{_9}(M)$& $0.754 $&$0.684 $& $0.573$& $0.525$& $0.436 $&$0.391 $& $0.342$& $0.286$ &$ 0.054 $\\
	\hline
	$r_{_{10}}(M)$& $0.755 $ &$0.684 $& $0.573$& $0.525$& $0.436 $&$0.391 $& $0.342$& $0.286$ &$ 0.054 $\\
	\hline
	\end{tabular}
	\vspace{1mm}
	\caption{The roots $ r_{_N}(M) $ for $ N=3,5, 6, 9, 10 $ $ M<{1}/{2(2\ln 2-1)} $.}
	\label{tabel-3}
\end{table}

 \begin{figure}[!htb]
 	\begin{center}
 		\includegraphics[width=0.50\linewidth]{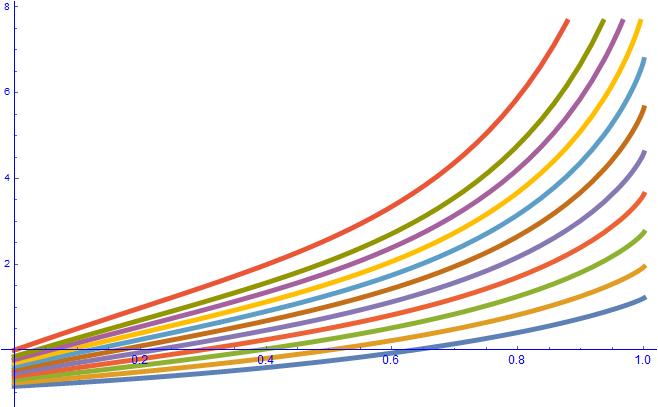}
 	\end{center}
 	\caption{The roots $ r_{_8}(M) $ of \eqref{e-2.6} for $ M<{1}/{2(2\ln 2-1)} $.}
 \end{figure}
We prove the following sharp Bohr-Rogosinski inequality by considering power of $ z $ in $ |f(z)| $ for the functions in the class $ \mathcal{P}^{0}_{\mathcal{H}}(M) $.
 \begin{thm}\label{th-2.7}
	Let $ f\in \mathcal{P}^{0}_{\mathcal{H}}(M) $ be given by \eqref{e-1.2}, then for $ N\geq 2 $,
	\begin{equation}\label{e-2.8}
		|f(z^m)|+\sum_{n=N}^{\infty}|a_n||z|^n\leq d\left(f(0),\partial f(\mathbb{D})\right)
	\end{equation} 
	for $ r\leq r_{_{m,N}}(M) $, where $ r_{_{m,N}}(M)\in (0,1) $ is the smallest root of the equation
	\begin{equation}\label{e-2.9}
		r^m-1-2M\left(r^m+r-1+\psi_m(r)+\psi_1(r)+\sum_{n=2}^{N-1}\frac{r^n}{n(n-1)}+\ln 4\right)=0,
	\end{equation}
	and $ \psi_m(r):=(1-r^m)\ln\,(1-r^m) $.	The radius $ r_{_{m,N}}(M) $ is the best possible.
\end{thm}
 
 \section{Improved Bohr radius for the class $ \mathcal{P}^{0}_{\mathcal{H}}(M) $}
 In $ 2017 $, Kayumov and Ponnusamy \cite{kayumov-2017} proved the following improved version of Bohr's inequality.
 \begin{thm}\cite{kayumov-2017}\label{th-1.3}
 	Let $f(z)=\sum_{n=0}^{\infty} a_{n}z^{n}$ be analytic in $\mathbb{D}$, $|f(z)|\leq 1$ and $S_{r}$ denote the image of the subdisk $|z|<r$ under mapping $f$. Then 
 	\begin{equation}
 		B_{1}(r):=\sum_{n=0}^{\infty}|a_n|r^n+ \frac{16}{9} \left(\frac{S_{r}}{\pi}\right) \leq 1 \quad \mbox{for} \quad r \leq \frac{1}{3}
 	\end{equation}
 	and the numbers $1/3$ and $16/9$ cannot be improved. Moreover, 
 	\begin{equation}
 		B_{2}(r):=|a_{0}|^{2}+\sum_{n=1}^{\infty}|a_n|r^n+ \frac{9}{8} \left(\frac{S_{r}}{\pi}\right) \leq 1 \quad \mbox{for} \quad r \leq \frac{1}{2}
 	\end{equation}
 	and the numbers $1/2$ and $9/8$ cannot be improved.
 \end{thm}
  Our aim is to prove a harmonic analogue of the Theorem \ref{th-1.3} for  functions $ f\in \mathcal{P}^{0}_{\mathcal{H}}(M)  $. It will be interesting to investigate Theorem \ref{th-1.3} in powers of $ {S_{r}}/{\pi} $. Therefore, in order to generalize Theorem \ref{th-1.3}, we consider a $ N^{\rm th} $ degree polynomial in $ S_r/\pi $ as follows 
 \begin{equation*}
 	P\left(\frac{S_{r}}{\pi}\right)=\left( \frac{S_{r}}{\pi}\right)^N+\left(\frac{S_{r}}{\pi}\right)^{N-1}+\cdots+\frac{S_{r}}{\pi}. 
 \end{equation*} 
 In this connection, we recall the polylogarithm function which is defined by a power series in $ z $, a Dirichlet series in $ s $.
 \begin{equation*}
 	Li_s(z)=\sum_{n=1}^{\infty}\frac{z^n}{n^s}=z+\frac{z^2}{2^s}+\frac{z^3}{3^s}+\cdots,
 \end{equation*} and this definition is valid for arbitrary complex order $ s $ and for all complex arguments $ z $ with $ |z|<1. $ Therefore, the dilogarithm, denoted as $ Li_2(z) $, is a particular case of the polylogarithm. 
 
 \begin{thm}\label{th-2.10}
 	Let $ f\in \mathcal{P}^{0}_{\mathcal{H}}(M) $ be given by \eqref{e-1.2}. Then 
 	\begin{equation}\label{e-2.14}
 		r+\sum_{n=2}^{\infty}\bigg(|a_n|+|b_n|\bigg)r^n+P\left(\frac{S_{r}}{\pi}\right)\leq d\left(f(0),\partial f(\mathbb{D})\right)
 	\end{equation} 
 	for $ r\leq r_{_N}(M) $, where $ P(w)=w^N+w^{N-1}+\cdots+w $, a polynomial in $ w $ of degree $ N-1 $, and $ r_{_N}(M)\in (0,1) $ is the smallest root of the equation
 	\begin{equation}\label{e-2.15}
 		r-1+2M\left(r-1+(1-r)\ln(1-r)+2\ln 2+P\left(r^2+4M^2G(r)\right)\right)=0,
 	\end{equation}
 	where $ G(r) $ is defined by
 	\begin{equation*}
 		G(r):=r^2(Li_2(r^2)-1)+(1-r^2)\ln\;(1-r^2).
 	\end{equation*}
 	The radius $ r_{_N}(M) $ is the best possible.
 \end{thm}
By considering powers in the coefficients, we prove the following sharp result for functions in the class $ \mathcal{P}^{0}_{\mathcal{H}}(M) $. 
 \begin{thm}\label{th-2.13}
 	Let $ f\in \mathcal{P}^{0}_{\mathcal{H}}(M) $ be given by \eqref{e-1.2}, then
 	\begin{equation}\label{e-2.11}
 		r+\sum_{n=2}^{\infty}\bigg(|a_n|+|b_n|+n(n-1)(|a_n|+|b_n|)^2\bigg)r^n\leq d\left(f(0),\partial f(\mathbb{D})\right),\; \text{for}\;\;  r\leq r_{_M},
 	\end{equation} 
 	where $ r_{_M}\in (0,1) $ is the smallest root of the equation
 	\begin{equation}\label{e-2.12}
 		r-1-2M(2M+1)\left(r+(1-r)\ln\,(1-r)\right)-2M(1-\ln 4)=0.
 	\end{equation}
 	The radius $ r_{_M} $ is the best possible.
 \end{thm}
 
 \section{Refined Bohr radius  and Bohr-type inequalities for the class $ \mathcal{P}^{0}_{\mathcal{H}}(M) $}
 In this section, we prove sharp results refined-Bohr radius and certain Bohr-type inequalities for the class $ \mathcal{P}^{0}_{\mathcal{H}}(M) $.
 \begin{thm}\label{th-2.19}
 	Let $ f\in \mathcal{P}^{0}_{\mathcal{H}}(M) $ be given by \eqref{e-1.2}. Then for integer $ p\geq 1 $, $ N\geq 2 $ and $ t=[(N-1)/2] $, we have
 	\begin{align}
 		\label{e-2.20} &
 		|f(z)|^p+\sum_{n=N}^{\infty}(|a_n|+|b_n|)r^n+sgn(t)\sum_{n=2}^{t}(|a_n|+|b_n|)^2\frac{r^N}{1-r}\\&\nonumber\quad\quad+\frac{1}{1-r}\sum_{n=t+1}^{\infty}n(n-1)(|a_n|+|b_n|)^2r^{2n}\\&\nonumber\leq d\left(f(0),\partial f(\mathbb{D})\right),\; \text{for}\;\;  r\leq r_{p,t,N}{M},
 	\end{align}
 	where $ r_{p,t,N}{M}\in (0,1) $ is the smallest root of the equation
 	\begin{align}&\label{e-2.21}
 		(r+2M(r+(1-r)\ln(1-r)))^p+2M\left(r+(1-r)\ln(1-r)-\sum_{n=2}^{N-1}\frac{r^n}{n(n-1)}\right)\\&\nonumber\quad+sgn(t)\sum_{n=2}^{t}\frac{2M}{n(n-1)}\frac{r^N}{1-r}+\frac{4M^2}{1-r}\left(r^2+(1-r^2)\ln(1-r^2)-\sum_{n =1}^{t}\frac{r^{2n}}{n(n-1)}\right)\\&\nonumber \quad-1-2M(1-\ln 4)=0.
 	\end{align}
 	The radius $ r_{p,t,N}(M) $ is the best possible.
 \end{thm}
Liu, Shang and Xu \cite{Liu-Shang-Xu-JIA-2018} have proved the following result computing Bohr-type radius for the analytic functions $ f(z) $ for which $ |a_0| $ and $ |a_1| $ are replaced by $ |f(z)| $ and $ f^{\prime}(z) $ respectively.
\begin{thm}\cite{Liu-Shang-Xu-JIA-2018}\label{th-5.1}
	Suppose that $ f(z)=\sum_{n=0}^{\infty}a_nz^n $ is analytic in $ \mathbb{D} $ and $ |f(z)|<1 $ in $ \mathbb{D} $. Then
	\begin{equation*}
		|f(z)|+|f^{\prime}(z)||z|+\sum_{n=2}^{\infty}|a_n||z|^n\leq\quad\mbox{for}\quad |z|=r\leq\frac{\sqrt{17}-3}{4}.
	\end{equation*}
The radius $ (\sqrt{17}-3)/4 $ is the best possible.
\end{thm} 
\noindent We prove the following sharp Bohr-type inequality by considering Jacobian in place of derivative for the functions of the class $ \mathcal{P}^{0}_{\mathcal{H}}(M) $ which is a harmonic analogue of Theorem \ref{th-5.1}.
 \begin{thm}\label{th-2.16}
 	For the function $ f\in \mathcal{P}^{0}_{\mathcal{H}}(M) $ given by \eqref{e-1.2}, we have
 	\begin{align}
 		\label{e-2.17} &
 		|f(z)|^p+\sqrt{{|J}_f(z)}|r+\sum_{n=2}^{\infty}(|a_n|+|b_n|)r^n+\frac{1}{1-r^N}\sum_{n=2}^{\infty}n(n-1)(|a_n|+|b_n|)^2r^{2n}\\&\nonumber\leq d\left(f(0),\partial f(\mathbb{D})\right),\; \text{for}\;\;  r\leq r_{p,N}(M),
 	\end{align}
 	where $ r_{p,N}(M)\in (0,1) $ is the smallest root of the equation
\begin{align}&\label{e-2.18}
	(r+2M(r+(1-r)\ln(1-r)))^p+(1-2M\ln(1-r))r+2M(r+(1-r)\ln(1-r))\\&\nonumber\quad\quad+\frac{4M^2}{1-r^N}(r^2+(1-r^2)\ln(1-r^2))-1-2M(1-\ln 4)=0.
\end{align}
 	The radius $ r_{p,N}(M) $ is the best possible.
 \end{thm}

\section{Proof of the main results}	
\noindent Before starting the proof of the main results, we recall here growth formula and distance bound for the class $ \mathcal{P}^{0}_{\mathcal{H}}(M) $. For $ f \in \mathcal{P}^{0}_{\mathcal{H}}(M) $, we have
\begin{equation}\label{e-3.1}
	|f(z)|\geq |z|+2M\sum\limits_{n=2}^{\infty}  \dfrac{2(-1)^{n-1}}{n(n-1)}|z|^{n} \quad \mbox{for } \quad |z|<1.
\end{equation}
Then the Euclidean distance between $f(0)$ and the boundary of $f(\mathbb{D})$ is given by 
\begin{equation}\label{e-3.2}
	d(f(0), \partial f(\mathbb{D}))= \liminf\limits_{|z|\rightarrow 1} |f(z)-f(0)|.
\end{equation}
Since $f(0)=0$, from \eqref{e-1.4} and \eqref{e-3.2} we obtain 
\begin{equation}\label{e-3.3}
	d(f(0), \partial f(\mathbb{D})) \geq 1+2M\sum\limits_{n=2}^{\infty}  \dfrac{(-1)^{n-1}}{n(n-1)}.
\end{equation}
\begin{pf}[\bf Proof of Theorem \ref{th-2.1}]
	Let $f \in \mathcal{P}^{0}_{\mathcal{H}}(M)$ be given by \eqref{e-1.2}. Using Lemmas \ref{lem-1.2} and \ref{lem-1.3}, for $ |z|=r_N(M) $, we obtain 
	\begin{align}\label{e-3.4}&
		|f(z)|+\sum_{n=N}^{\infty}(|a_n|+|b_n|)r^n\\ &\leq\nonumber r+\sum_{n=2}^{\infty}\frac{2Mr^n}{n(n-1)}+\sum_{n=N}^{\infty}\frac{2Mr^n}{n(n-1)}\\&=\nonumber r+2M(r+(1-r)\ln(1-r))+2M\left(r+(1-r)\ln(1-r)-\sum_{n=2}^{N-1}\frac{r^n}{n(n-1)}\right)\\&=\nonumber r+2M\left(2r+2(1-r)\ln(1-r)-\sum_{n=2}^{N-1}\frac{r^n}{n}\right).
	\end{align}
	A simple computation shows that
	\begin{equation}\label{e-3.5}
		r+2M\left(2r+2(1-r)\ln(1-r)-\sum_{n=2}^{N-1}\frac{r^n}{n(n-1)}\right)\leq 1+2M(1-2\ln 2)
	\end{equation}
	for $ r\leq r_{_N}(M) $, where $ r_{_N}(M) $ is the smallest root of $ F_1(r)=0 $ in $ (0,1) $, here
$ F_1 : [0,1)\rightarrow \mathbb{R} $ is defined by 
	$$ F_1(r):=r-1+2M\left(2r-1+2(1-r)\ln(1-r)-\sum_{n=2}^{N-1}\frac{r^n}{n(n-1)}+2\ln 2\right).
	$$
	\par The existence of the root $ r_{_N}(M) $ is guaranteed by the following fact that $ F_1 $ is a continuous function with the properties $ F_1(0)=-1-2M(1-2\ln 2)<0 $ and $ \displaystyle\lim_{r\rightarrow 1}F_1(r)=+\infty. $ Let $ r_{_N}(M) $ to be the smallest root of $ F_1(r)=0 $ in $ (0,1) $. Then we have $ F_1(r_{_N}(M))=0 $. That is
	\begin{align}\label{e-3.6} 
		r_{_N}(M)-1+2MJ_{_{N,M}}(r)=0,
	\end{align}
	where
	\begin{equation*}
		J_{_{N,M}}(r)=2r_{_N}(M)-1+2(1-r_{_N}(M))\ln(1-r_{_N}(M))-\sum_{n=2}^{N-1}\frac{r^n_{_N}(M)}{n(n-1)}+2\ln 2.
	\end{equation*}
\noindent It follows from \eqref{e-3.3}, \eqref{e-3.4} and \eqref{e-3.5} for $ |z|=r\leq r_{_N}(M) $, we obtain
	$$ 
	|f(z)|+\sum_{n=N}^{\infty}(|a_n|+|b_n|)r^n\leq d(f(0),\partial f(\mathbb{D})).
	$$ 
	In order to show that $ r_{_N}(M) $ is the best possible radius, we consider the following function $ f=f_{M} $ be defined by 
	\begin{equation}\label{e-3.7a}
		f_{M}(z)=z+\sum_{n=2}^{\infty}\frac{2Mz^n}{n(n-1)}.
	\end{equation}
	It is easy to show that $ f_{M}\in\mathcal{P}^{0}_{\mathcal{H}}(M) $ and for $ f=f_{M} $, we have 
	\begin{equation}\label{e-3.7b}
		d(f(0),\partial f(\mathbb{D}))=1+2M(1-2\ln 2).
	\end{equation}
	For the function $ f=f_{M} $ and $ |z|=r_{_N}(M) $, a simple computation using \eqref{e-3.5} and \eqref{e-3.7a} shows that 
	\begin{align*} &
		|f(z)|+\sum_{n=N}^{\infty}(|a_n|+|b_n|)r^n_{_N}(M)\\ &= r_{_N}(M)+2M(r_N(M)+(1-r_{_N}(M))\ln(1-r_{_N}(M)))\\&\quad\quad+2M\left(r_{_N}(M)+(1-r_{_N}(M))\ln(1-r_{_N}(M))-\sum_{n=2}^{N-1}\frac{r^n_{_N}(M)}{n(n-1)}\right)\\&=\nonumber r_{_N}(M)+2M\left(2r_{_N}(M)+2(1-r_{_N}(M))\ln(1-r_{_N}(M))-\sum_{n=2}^{N-1}\frac{r^n_{_N}(M)}{n}\right)
		\\&\nonumber= 1+2M(1-2\ln 2)\\&=d(f(0),\partial f(\mathbb{D})).
	\end{align*}
	Therefore, the radius $ r_{_N}(M) $ is the best possible. This completes the proof.
\end{pf}

\begin{pf}[\bf Proof of Theorem \ref{th-2.4}]
	Let $f \in \mathcal{P}^{0}_{\mathcal{H}}(M)$ be given by \eqref{e-1.2}. 
	Using Lemmas \ref{lem-1.2} and \ref{lem-1.3} for $ |z|=r $, we obtain 
	\begin{align}\label{e-3.8}&
		|f(z)|^2+\sum_{n=N}^{\infty}(|a_n|+|b_n|)r^n\\ &\leq\nonumber \left(r+\sum_{n=2}^{\infty}\frac{2Mr^n}{n(n-1)}\right)^2+\sum_{n=N}^{\infty}\frac{2Mr^n}{n(n-1)}\\&=\nonumber (r+2M(r+(1-r)\ln(1-r)))^2+2M\left(r+(1-r)\ln(1-r)-\sum_{n=2}^{N-1}\frac{r^n}{n(n-1)}\right).
	\end{align}
	An elementary calculation shows that
\begin{align}\label{e-3.9} &
	(r+2M(r+(1-r)\ln(1-r)))^2+2M\left(r+(1-r)\ln(1-r)-\sum_{n=2}^{N-1}\frac{r^n}{n(n-1)}\right)\\&\nonumber\leq 1+2M(1-2\ln 2)
\end{align}
	for $ r\leq r_{_N}(M) $, where $ r_{_N}(M) $ is the smallest root of $ F_2(r)=0 $ in $ (0,1) $, here
	\begin{align*} 
		F_2(r):&=(r+2M(r+(1-r)\ln(1-r)))^2-1\\&\quad\quad+2M\left(r-1+(1-r)\ln(1-r)-\sum_{n=2}^{N-1}\frac{r^n}{n(n-1)}+2\ln 2\right).
	\end{align*}
 Therefore, we have $ F_2(r_{_N}(M))=0 $. That is 
	\begin{align}\label{e-3.10} &
		\bigg(r_{_N}(M)+2M(r_{_N}(M)+(1-r_{_N}(M))\ln(1-r_{_N}(M)))\bigg)^2-1\\&\nonumber\quad\quad+2M\left(r_{_N}(M)-1+(1-r_{_N}(M))\ln(1-r_{_N}(M))-\sum_{n=2}^{N-1}\frac{r^n_{_N}(M)}{n(n-1)}+2\ln 2\right)=0.
	\end{align}
\noindent From \eqref{e-3.3}, \eqref{e-3.8} and \eqref{e-3.9} for $ |z|=r\leq r_{_N}(M) $, we obtain
	$$ 
	|f(z)|^2+\sum_{n=N}^{\infty}(|a_n|+|b_n|)r^n\leq d(f(0),\partial f(\mathbb{D})).
	$$ 
	To show that the radius $ r_{_N}(M) $ is the best possible, we consider the function $ f=f_{M} $ given by \eqref{e-3.7a}. For the function $ f=f_{M} $ and $ |z|=r_{_N}(M) $, a simple computation  using \eqref{e-3.10} and \eqref{e-3.7b} shows that 
	\begin{align*} &
		|f(z)|^2+\sum_{n=N}^{\infty}(|a_n|+|b_n|)r^n_{_N}(M)\\ &= 	(r_{_N}(M)+2M(r_{_N}(M)+(1-r_{_N}(M))\ln(1-r_{_N}(M))))^2\\&\quad\quad+2M\left(r_{_N}(M)+(1-r_{_N}(M))\ln(1-r_{_N}(M))-\sum_{n=2}^{N-1}\frac{r^n_{_N}(M)}{n(n-1)}\right)\\&\nonumber= 1+2M(1-2\ln 2)\\&=d(f(0),\partial f(\mathbb{D})).
	\end{align*}
	Hence, the radius $ r_{_N}(M) $ is the best possible. This completes the proof.
\end{pf}	

\begin{pf}[\bf Proof of Theorem \ref{th-2.7}]
	Let $f \in \mathcal{P}^{0}_{\mathcal{H}}(M)$ be given by \eqref{e-1.2}. 	Using Lemmas \ref{lem-1.2} and \ref{lem-1.3} for $ |z|=r $, we obtain 
\begin{align}\label{e-3.12}&
|f(z^m)|+\sum_{n=N}^{\infty}|a_n||z|^n\\ &\leq\nonumber r^m+\sum_{n=2}^{\infty}\frac{2M(r^m)^n}{n(n-1)}+\sum_{n=N}^{\infty}\frac{2Mr^n}{n(n-1)}\\&=\nonumber r^m+2M(r^m+(1-r^m)\ln(1-r^m))+2M\left(r+(1-r)\ln(1-r)-\sum_{n=2}^{N-1}\frac{r^n}{n(n-1)}\right).
\end{align}
	A simple computation shows that
\begin{align}
	\label{e-3.13} & r^m+2M(r^m+(1-r^m)\ln(1-r^m))+2M\left(r+(1-r)\ln(1-r)-\sum_{n=2}^{N-1}\frac{r^n}{n(n-1)}\right)\\&\nonumber\leq 1+2M(1-2\ln 2)
\end{align}
	for $ r\leq r_{_N}(M) $, where $ r_{_N}(M) $ is the smallest root of $ F_3(r)=0 $ in $ (0,1) $, here 
	\begin{align*} 
		F_3(r):&=r^m-1+2M\left(r^m+r-1+(1-r^m)\ln(1-r^m)+(1-r)\ln(1-r)\right)\\&\quad\quad-2M\left(\sum_{n=2}^{N-1}\frac{r^n}{n(n-1)}-2\ln 2\right).
	\end{align*}
	 Therefore, we have $ F_3(r_{_{m,N}}(M))=0 $. That is 
	\begin{align}\label{e-3.14} 
		&2M\left(r^m_{_{m,N}}(M)+r-1+(1-r^m_{_{m,N}}(M))\ln(1-r^m_{_{m,N}}(M))+(1-r_{_{m,N}}(M))\ln(1-r_{_{m,N}}(M))\right)\\&\nonumber \quad\quad+ r^m_{_{m,N}}(M)-1-2M\left(\sum_{n=2}^{N-1}\frac{r^n_{_{m,N}}(M)}{n(n-1)}-2\ln 2\right)=0.
	\end{align}
	From \eqref{e-3.3}, \eqref{e-3.12} and \eqref{e-3.13} for $ |z|=r\leq r_{_{m,N}}(M) $, we obtain
	$$ 
	|f(z^m)|+\sum_{n=N}^{\infty}(|a_n|+|b_n|)r^n\leq d(f(0),\partial f(\mathbb{D})).
	$$ 
	To show that $ r_{_{m,N}}(M) $ is the best possible radius, we consider the following function $ f=f_{M} $ defined by \eqref{e-3.7a}. For $ f=f_{M} $ and $ |z|=r_{_{m,N}}(M) $, a simple computation using \eqref{e-3.14} and \eqref{e-3.7b} shows that 
	\begin{align*} &
		|f(z)|+\sum_{n=N}^{\infty}(|a_n|+|b_n|)r^n_{_N}(M)\\ &= r_{_{m,N}}(M)+2M(r_N(M)+(1-r_{_{m,N}}(M))\ln(1-r_{_{m,N}}(M)))\\&\quad\quad+2M\left(r_{_{m,N}}(M)+(1-r_{_{m,N}}(M))\ln(1-r_{_{m,N}}(M))-\sum_{n=2}^{N-1}\frac{r^n_{_{m,N}}(M)}{n(n-1)}\right)\\&=\nonumber r_{_{m,N}}(M)+2M\left(2r_{_{m,N}}(M)+2(1-r_{_{m,N}}(M))\ln(1-r_{_{m,N}}(M))-\sum_{n=2}^{N-1}\frac{r^n_{_{m,N}}(M)}{n(n-1)}\right)
		\\&\nonumber= 1+2M(1-2\ln 2)\\&=d(f(0),\partial f(\mathbb{D})).
	\end{align*}
	Hence, the radius $ r_{_{m,N}}(M) $ is the best possible. This completes the proof.
\end{pf}	

\begin{pf}[\bf Proof of Theorem \ref{th-2.10}]
	Let $f \in \mathcal{P}^{0}_{\mathcal{H}}(M)$ be given by \eqref{e-1.2}. 
 It is well known that 
 \begin{align}\label{e-5.15}
 	S_r=\iint\limits_{D_r}&\left(|h^{\prime}(z)|^2-|g^{\prime}(z)|^2\right)dxdy,\\ \frac{1}{\pi}\iint\limits_{D_r}& \label{e-5.16}|h^{\prime}(z)|^2dxdy=\sum_{n=1}^{\infty}n|a_n|^2r^{2n}\\ \frac{1}{\pi}\iint\limits_{D_r}&\label{e-5.17}|g^{\prime}(z)|^2dxdy=\sum_{n=2}^{\infty}n|b_n|^2r^{2n}.
 \end{align}
 Then, we in view of Lemma \ref{lem-1.2} and using \eqref{e-5.15}, \eqref{e-5.16} and \eqref{e-5.17}, we obtain
 \begin{align*}
 	\frac{S_r}{\pi}&=\frac{1}{\pi}\iint\limits_{D_r}\left(|h^{\prime}(z)|^2-|g^{\prime}(z)|^2\right)dxdy\\&=r^2+\sum_{n=2}^{\infty}|a_n|^2r^{2n}-\sum_{n=2}^{\infty}n|b_n|^2r^{2n}\\&= r^2+\sum_{n=2}^{\infty}n\left(|a_n|+|b_n|\right)\left(|a_n|-|b_n|\right)r^{2n}\\&\leq r^2+\sum_{n=2}^{\infty}\frac{4M^2r^{2n}}{n^2(n-1)^2}\\&=r^2+4M^2\left(r^2(Li_2(r^2)-1)+(r^2-1)\ln (1-r^2)\right).
 \end{align*}
	Using Lemmas \ref{lem-1.2} and \ref{lem-1.3} for $ |z|=r $, we obtain 
	\begin{align}\label{e-3.16}
	r+\sum_{n=2}^{\infty}(|a_n|+|b_n|)r^n+P\left(\frac{S_{r}}{\pi}\right)&\leq r+\sum_{n=2}^{\infty}\frac{2Mr^n}{n(n-1)}+P(r^2+4M^2G_1(r))\\&=\nonumber r+2M(r+(1-r)\ln(1-r))+P(r^2+4M^2G(r)).
	\end{align}
where $ G_1(r) $ is defined by
\begin{equation*}
	G_{1}(r):=r^2(Li_2(r^2)-1)+(1-r^2)\ln\;(1-r^2).
\end{equation*}
	A simple computation shows that
	\begin{equation}\label{e-3.17}
		 r+2M(r+(1-r)\ln(1-r))+P(r^2+4M^2G_{1}(r))\leq 1+2M(1-2\ln 2)
	\end{equation}
	for $ r\leq r_{_N}(M) $, where $ r_{_N}(M) $ is the smallest root of $ F_4(r)=0 $ in $ (0,1) $, here 
	$$ F_4(r):=	r-1+2M(r-1+(1-r)\ln(1-r)+P(r^2+4M^2G_1(r))+2\ln 2).
	$$
	Therefore, we have $ F_4(r_{_N}(M))=0 $. Then we have 
\begin{align}
	\label{e-3.18} & r(M)-1+2M\left(r(M)-1+(1-r(M))\ln(1-r(M))\right)\\&\nonumber\quad\quad+2M\left(P(r^2(M)+4M^2G_1(r(M)))+2\ln 2\right)=0.
\end{align}
	
	\noindent Using \eqref{e-3.3}, \eqref{e-3.17} and \eqref{e-3.18} for $ |z|=r\leq r_{_N}(M) $, we obtain
	$$ 
	|f(z)|+\sum_{n=N}^{\infty}(|a_n|+|b_n|)r^n\leq d(f(0),\partial f(\mathbb{D})).
	$$ 
	To show that $ r_{_N}(M) $ is the best possible radius, we consider the function $ f=f_{M} $ defined by \eqref{e-3.7a}. 
 Then it is easy to see that
	\begin{equation}\label{e-3.19}
		d(f(0),\partial f(\mathbb{D}))=1+2M(1-2\ln 2).
	\end{equation}
	For $ f=f_{M} $ and $ |z|=r_{_N}(M) $, a simple computation using \eqref{e-3.18} and \eqref{e-3.19} shows that 
	\begin{align*} &
		|f(z)|+\sum_{n=N}^{\infty}(|a_n|+|b_n|)r^n(M)\\ &= r(M)+2M(r_{_N}(M)+(1-r(M))\ln(1-r_{_N}(M)))\\&\quad\quad+2M\left(r(M)+(1-r(M))\ln(1-r_{_N}(M))-\sum_{n=2}^{N-1}\frac{r^n_{_N}(M)}{n(n-1)}\right)\\&=\nonumber r(M)+2M\left(2r_{_N}(M)+2(1-r_{_N}(M))\ln(1-r_{_N}(M))-\sum_{n=2}^{N-1}\frac{r^n_{_N}(M)}{n}\right)
		\\&\nonumber= 1+2M(1-2\ln 2)\\&=d(f(0),\partial f(\mathbb{D})).
	\end{align*}
	Hence the radius $ r_{_N}(M) $ is the best possible. This completes the proof.
\end{pf}	

\begin{pf}[\bf Proof of Theorem \ref{th-2.13}]
	Let $f \in \mathcal{P}^{0}_{\mathcal{H}}(M)$ be given by \eqref{e-3.7a}.	Using Lemmas \ref{lem-1.2} and \ref{lem-1.3} for $ |z|=r $, we obtain 
	\begin{align}\label{e-3.20}
		r+\sum_{n=2}^{\infty}&\left(|a_n|+|b_n|+n(n-1)(|a_n|+|b_n|)^2\right)r^n\\ &\leq\nonumber r+\sum_{n=2}^{\infty}\frac{2Mr^n}{n(n-1)}+\sum_{n=2}^{\infty}\frac{4M^2r^n}{n(n-1)}\\&=\nonumber r-2M(1+2M)(r+(1-r)\ln(1-r)).
	\end{align}
	An simple calculation shows that
	\begin{equation}\label{e-3.21}
		r-2M(1+2M)(r+(1-r)\ln(1-r))\leq 1+2M(1-2\ln 2)
	\end{equation}
	for $ r\leq r_{_M} $, where $ r_{_M} $ is the smallest root of $ F_5(r)=0 $ in $ (0,1) $, here
	$$ F_5(r):=r-1-2M(2m+1)(r+(1-r)\ln (1-r))-2M(1-2\ln 2).
	$$
	 Then we have $ F_5(r_{_M})=0 $. That is,
	\begin{align}\label{e-3.22} &
		r_{_M}-1-2M(2m+1)(r_{_M}+(1-r_{_M})\ln (1-r_{_M}))-2M(1-2\ln 2)=0.
	\end{align}
	
	\noindent From \eqref{e-3.3}, \eqref{e-3.20} and \eqref{e-3.21} for $ |z|=r\leq r_{_M} $, we obtain
	$$ 
	r+\sum_{n=2}^{\infty}\left(|a_n|+|b_n|+n(n-1)(|a_n|+|b_n|)^2\right)r^n\leq d(f(0),\partial f(\mathbb{D})).
	$$ 
	In order to show that $ r_{_M} $ is the best possible radius, we consider the function $ f=f_{M} $ defined by \eqref{e-3.7a}. For the function $ f=f_{M} $ and $ |z|=r_{_M} $, a simple computation using \eqref{e-3.22} and \eqref{e-3.7b} shows that 
	\begin{align*} &
			r_{_M}+\sum_{n=2}^{\infty}\left(|a_n|+|b_n|+n(n-1)(|a_n|+|b_n|)^2\right)r^n_{_M}\\ &= r_{_M}-2M(1+2M)(r_{_M}+(1-r_{_M})\ln(1-r_{_M}))\\&\nonumber= 1+2M(1-2\ln 2)\\&=d(f(0),\partial f(\mathbb{D})).
	\end{align*}
	This shows that the radius $ r_{_M} $ is the best possible. This completes the proof.
\end{pf}	

\begin{pf}[\bf Proof of Theorem \ref{th-2.16}]
	Let $f \in \mathcal{P}^{0}_{\mathcal{H}}(M)$ be given by \eqref{e-1.2}. The Jacobian of complex-valued harmonic function $ f=h+\overline{g} $ has the following property 
\begin{align}\label{e-5.25}
	|{J}_f(z)|\leq |h^{\prime}(z)|^2-|g^{\prime}(z)|^2\leq |h^{\prime}(z)|^2\leq \left(1+2M\sum_{n=1}^{\infty}\frac{r^n}{n}\right)^2.
\end{align}
	Using Lemmas \ref{lem-1.2} and \ref{lem-1.3} and \eqref{e-5.25} for $ |z|=r $, we obtain 
\begin{align}\label{e-3.24}&
	|f(z)|^p+r\sqrt{{|J}_f(z)}|+\sum_{n=2}^{\infty}(|a_n|+|b_n|)r^n+\frac{1}{1-r^N}\sum_{n=2}^{\infty}n(n-1)(|a_n|+|b_n|)^2r^{2n}\\ &\leq\nonumber \left(r+\sum_{n=2}^{\infty}\frac{2Mr^n}{n(n-1)}\right)^p+r\left(1+2M\sum_{n=1}^{\infty}\frac{(r^2)^n}{n}\right)+\sum_{n=2}^{\infty}\frac{2Mr^n}{n(n-1)}\\&\nonumber\quad\quad+\frac{1}{1-r^N}\sum_{n=2}^{\infty}\frac{4M^2(r^2)^n}{n(n-1)}\\&=\nonumber (r+2M(r+(1-r)\ln(1-r)))^p+(1-2M(\ln(1-r)))r\\&\nonumber\quad\quad+2M\left(r+(1-r)\ln(1-r)\right)+\frac{4M^2}{1-r^N}\left(r^2+(1-r^2)\ln(1-r^2)\right)\\&:=B_{M,N,p}(r)\nonumber.
\end{align}
A simple computation shows that
\begin{align}
	\label{e-3.25}
	B_{M,N,p}(r)\leq 1+2M(1-2\ln 2)
\end{align}	for $ r\leq r_{_{p,N}}(M) $, where $ r_{_{p,N}}(M) $ is the smallest root of $ F_6(r)=0 $ in $ (0,1) $, here 
	$$ F_6(r):=B_{M,N}(r)- 1-2M(1-2\ln 2).
	$$ 
	Therefore, we have $ F_6(r_{_{p,N}}(M))=0 $. That is 
	\begin{align}\label{e-3.28} &
	B_{M,N}(r_{_{p,N}}(M))- 1-2M(1-2\ln 2)=0.
	\end{align}
	
	\noindent From \eqref{e-3.3}, \eqref{e-3.24} and \eqref{e-3.25} for $ |z|=r\leq r_{_{p,N}}(M) $, we obtain
	\begin{align*} &
		|f(z)|^p+r\sqrt{{|J}_f(z)}|+\sum_{n=2}^{\infty}(|a_n|+|b_n|)r^n+\frac{1}{1-r^N}\sum_{n=2}^{\infty}n(n-1)(|a_n|+|b_n|)^2r^{2n}\\&\leq d(f(0),\partial f(\mathbb{D})).
	\end{align*} 
	To show that $ r_{_{p,N}}(M) $ is the best possible radius, we consider the function $ f=f_{M} $ defined by \eqref{e-3.7a}.	For $ f=f_{M} $ and $ |z|=r_{_{p,N}}(M) $, a simple calculation using \eqref{e-3.28} and \eqref{e-3.7b} shows that 
	\begin{align*} &
	|f(z)|^p+r_{_{p,N}}(M)\sqrt{{|J}_f(z)}|+\sum_{n=2}^{\infty}(|a_n|+|b_n|)r^n_{_{p,N}}(M)\\&\quad\quad+\frac{1}{1-r^N_{_{p,N}}(M)}\sum_{n=2}^{\infty}n(n-1)(|a_n|+|b_n|)^2r^{2n}_{_{p,N}}(M)
		\\&= 	B_{M,N,p}(r_{_{p,N}}(M))\nonumber= 1+2M(1-2\ln 2)\\&=d(f(0),\partial f(\mathbb{D})).
	\end{align*}
	Therefore, the radius $ r_{_{p,N}}(M) $ is the best possible. This completes the proof.
\end{pf}	

\begin{pf}[\bf Proof of Theorem \ref{th-2.19}]
	Let $f \in \mathcal{P}^{0}_{\mathcal{H}}(M)$ be given by \eqref{e-1.2}. 
	Using Lemmas \ref{lem-1.2} and \ref{lem-1.3} for $ |z|=r $, we obtain 
	\begin{align}\label{e-3.30}&
	|f(z)|^p+\sum_{n=N}^{\infty}(|a_n|+|b_n|)r^n+sgn(t)\sum_{n=2}^{t}(|a_n|+|b_n|)^2\frac{r^N}{1-r}\\&\nonumber\quad\quad+\frac{1}{1-r}\sum_{n=t+1}^{\infty}n(n-1)(|a_n|+|b_n|)^2r^{2n}\\ &\leq\nonumber \left(r+2M(r+(1-r)\ln(1-r))\right)^p+2M\left(r+(1-r)\ln(1-r)-\sum_{n=2}^{N-1}\frac{r^n}{n(n-1)}\right)\\&\nonumber\quad\quad+sgn(t)\sum_{n=2}^{t}\frac{2M}{n(n-1)}\frac{r^N}{1-r}+\frac{4M^2}{1-r}\left(r^2+(1-r^2)\ln(1-r^2)-\sum_{n=2}^{t}\frac{r^{2n}}{n(n-1)}\right)\\:&= C_{p,t,N,M}(r).\nonumber
	\end{align}
	It is easy to see that 
	\begin{equation}\label{e-3.31}
		C_{p,t,N,M}(r)\leq 1+2M(1-2\ln 2)
	\end{equation}
	for $ r\leq r_{_{p,t,N}}(M) $, where $ r_{_{p,t,N}}(M) $ is the smallest root of $ F_7(r)=0 $ in $ (0,1) $, here $ F_7 : [0,1)\rightarrow \mathbb{R} $ is defined by 
	$$ F_7(r):=C_{p,t,N,M}(r)- 1-2M(1-2\ln 2).
	$$
	Clearly, $ F_7(r_{_{p,t,N}}(M))=0 $. That is 
	\begin{align}\label{e-3.32} &
	C_{p,t,N,M}(r_{_{p,t,N}}(M))- 1-2M(1-2\ln 2)=0.
	\end{align}
	
	\noindent From \eqref{e-3.3}, \eqref{e-3.30} and \eqref{e-3.31} for $ |z|=r\leq r_{_{p,t,N}}(M) $, we obtain
	$$ 
	|f(z)|+\sum_{n=N}^{\infty}(|a_n|+|b_n|)r^n\leq d(f(0),\partial f(\mathbb{D})).
	$$ 
	In order to show that $ r_{_{p,t,N}}(M) $ is the best possible radius, we consider the function $ f=f_{M} $ defined by \eqref{e-3.7a}. For $ f=f_{M} $ and $ |z|=r_{_{p,t,N}}(M) $, a simple computation using \eqref{e-3.32} and \eqref{e-3.7b} shows that 
	\begin{align*} &
	|f(z)|^p+\sum_{n=N}^{\infty}(|a_n|+|b_n|)r^n_{_{p,t,N}}(M)+sgn(t)\sum_{n=2}^{t}(|a_n|+|b_n|)^2\frac{r^N_{_{p,t,N}}(M)}{1-r_{_{p,t,N}}(M)}\\&\nonumber\quad\quad+\frac{1}{1-r_{_{p,t,N}}(M)}\sum_{n=t+1}^{\infty}n(n-1)(|a_n|+|b_n|)^2(r^{2n}_{_{p,t,N}}(M))\\&=C_{p,t,N,M}(r_{_{p,t,N}}(M))= 1+2M(1-2\ln 2)\\&=d(f(0),\partial f(\mathbb{D})).
	\end{align*}
	Hence the radius $ r_{_{p,t,N}}(M) $ is the best possible. This completes the proof.
\end{pf}

\noindent\textbf{Acknowledgment:}  The first author is supported by the Institute Post Doctoral Fellowship of IIT Bhubaneswar, India, the second author is supported by SERB-MATRICS, and third author is supported by CSIR, India.

\end{document}